\documentclass[a4paper,12pt]{amsart}
\usepackage{amssymb,amsmath,hyperref,eufrak,bbm}
\author[Florent Benaych-Georges]{Florent Benaych-Georges}\address{Florent Benaych-Georges, LPMA,  UPMC Univ Paris 6, Case courier 188, 4, Place Jussieu, 75252 Paris Cedex 05, France.} \email{florent.benaych@upmc.fr}
\title[Limit theorems for the  unitary Brownian motion]{Finite dimensional projections of the  Brownian motion on large unitary groups}

\thanks{This work was partially supported by the \emph{Agence Nationale de la
Recherche} grant ANR-08-BLAN-0311-03.}
\date{\today}
\newcommand{\bbm}{\begin{bmatrix}}
\newcommand{\ebm}{\end{bmatrix}}
\newcommand{\ltoLtwo}{\stackrel{L^2}{\lto}}
\newcommand{\lto}{\longrightarrow}

\newcommand{\ovl}{\overline}
\newcommand{\bes}{\begin{equation*}}
\newcommand{\ees}{\end{equation*}}
\newcommand{\be}{\begin{equation}}
\newcommand{\ee}{\end{equation}}
\newcommand{\beqy}{\begin{eqnarray}}
\newcommand{\eeqy}{\end{eqnarray}}
\newcommand{\beq}{\begin{eqnarray*}}
\newcommand{\eeq}{\end{eqnarray*}}

\newcommand{\lan}{\langle}
\newcommand{\ran}{\rangle}

\newcommand{\Tr}{\operatorname{Tr}}

\newcommand{\ninf}{\underset{n\to\infty}{\longrightarrow}}

\newcommand{\one}{\mathbbm{1}}

\newcommand{\E}{\mathbb{E}}

\newcommand{\C}{\mathbb{C}}

\newcommand{\ud}{\mathrm{d}}

\newcommand{\pro}{probability }

\newcommand{\f}{\frac}
\newcommand{\ff}{\frac{1}}
\newcommand{\lf}{\left}
\newcommand{\ri}{\right}

\newcommand{\st}{such that }
\newcommand{\la}{\lambda}

\newcommand{\ste}{\, ;\, }
\newcommand{\mc}{\mathcal }
\newcommand{\al}{\alpha}
\newcommand{\eps}{\varepsilon}

\newtheorem{Th}{Theorem}[section]

\newtheorem{lem}[Th]{Lemma}

\newtheorem{exs}[Th]{Examples}
\newtheorem{rmq}[Th]{Remark}
\newtheorem{cor}[Th]{Corollary}

\newenvironment{pr}{\noindent {\it Proof. }}{\hfill$\square$}

\long\def\symbolfootnote[#1]#2{\begingroup
\def\thefootnote{\fnsymbol{footnote}}\footnote[#1]{#2}\endgroup}

\begin{document}
\maketitle

\begin{abstract}In this paper, we are concerned with the large $n$ limit of the distributions of  linear combinations of the entries of a  Brownian motion on the group of $n\times n$ unitary matrices. We prove that the process of such a linear combination converges to a Gaussian one. Various scales of time and various initial distributions are concerned, giving rise to various limit processes,  related to the geometric construction of the unitary Brownian motion. As an application, we propose a quite  short proof of the asymptotic Gaussian feature of the linear combinations of the entries of Haar distributed random unitary matrices, a result already proved by Diaconis {\it et al}.\end{abstract}

\section*{Introduction}
There is a natural definition of Brownian motion on any compact Lie group, whose distribution is sometimes called the heat kernel measure. 
Mainly due to its relations with the object from free probability  theory called the free unitary Brownian motion and with the two-dimentional Yang-Mills theory, the Brownian motion on large unitary groups has appeared  in several papers during the last decade. Rains, in \cite{rains97}, Xu, in \cite{xu97}, Biane, in \cite{b97,b97b} and L\'evy and   Ma\"\i da, in \cite{SW,thierry-mylene}, are all concerned with the asymptotics of the spectral distribution of large random matrices distributed according to the heat kernel measure. Also, in \cite{demni2}, Demni makes use of the unitary Brownian motion in the study of Jacobi  processes, and, in \cite{benlevy08}, L\'evy and the author of the present paper construct a continuum of convolutions between the classical and free ones based on the conjugation of random matrices with a unitary Brownian motion. 
In this paper, we are concerned with the asymptotic distributions  of linear combinations of the entries of an $n\times n$ unitary Brownian motion as $n$  tends to infinity. 

We first give the joint limit distribution, as $n$ tends to infinity,  of (possibly rescaled) random processes of the type $(\Tr[A(V_t-I)])_{t\geq 0}$ for $(e^{-t/2}V_t)_{t\geq 0}$ a Brownian motion  staring at $I$ on the group of unitary $n\times n$ matrices and $A$ an   $n\times n$ matrix (Theorem \ref{17209.0543}). 
This theorem is the key result of the paper, since specifying the choice of the matrices $A$ and randomizing them will then allow us to prove all other results. As a first example, it allows us to find out, for any sequence $(\al_n)_n$ of positive numbers with a limit $\al\in [0,+\infty]$, the limit distribution of  any upper-left corner of $\sqrt{n/\al_n}(V_{\log(\al_nt+1)}-I)_{t\geq 0}$ (Corollary \ref{26209.16h29}):   for small scales of time (i.e. when $\al=0$), the limit process is purely skew-Hermitian, whereas for large scales of time ($\al=+\infty$), the limit process is a standard complex matricial Brownian motion. For intermediate scales of time ($0<\al<+\infty$), the limit process is an interpolation between these extreme cases. The existence of these three asymptotic regimes can be explained by the fact that 
the unitary Brownian motion is the ``wrapping", on the unitary group, of a Brownian motion on the tangent space of this group  at $I$ (which is the space of skew-Hermitian matrices), and that as the time goes to infinity, its distribution tends to the Haar measure (for which  the upper-left corners are asymptotically distributed as standard complex Gaussian random matrices).

Secondly, we consider a unitary Brownian motion $(e^{-t/2}V_t)_{t\geq 0}$ whose initial distribution is the uniform measure  on the group of permutation matrices: its rows are exchangeable, as its columns.    In this case, for any positive sequence $(\al_n)$ and any positive integer $p$,   the $p\times p$ upper left corner of $(\sqrt{n/\al_n}V_{\log(\al_nt+1)})_{t\geq 0}$ converges  to a standard complex matricial Brownian motion (Corollary \ref{26209.16h50}). 

Since  
the unitary Brownian motion  distributed according to the Haar measure at time zero has a stationary distribution,  our results allow us to give very short proofs of some well-known results of  Diaconis  {\it et al}, first proved in \cite{diaconis2003},  about the  asymptotic normality of linear combinations of the entries of uniform random unitary matrices (Theorem \ref{17209.0544}  and  Corollary \ref{27309.23h50}). 


It is clear that the same analysis would give similar results for the Brownian motion on the orthogonal group. For notational  brevity, we chose to focus on the unitary group.

Let us now present briefly  what problems underlie the asymptotics of linear combinations of the entries of a unitary Brownian motion.

{\it Asymptotic normality of random unit  vectors and unitary matrices:} The historical first result in this direction is due to  
\'Emile Borel, who proved a century ago, in  \cite{b1906}, that, for 
a uniformly distributed point $ (X_1 , \ldots , X_n)$ on the unit euclidian sphere $\mathbb{S}^{n-1}$, the scaled first coordinate $\sqrt{n}X_1$ converges weakly to the standard Gaussian distribution as the dimension $n$ tends to infinity. As  explained in the introduction of the paper \cite{diaconis2003}  of  Diaconis  {\it et al.},    this says that the features of the ``microcanonical" ensemble in a certain model for statistical mecanics (uniform measure on the sphere) are captured by the ``canonical" ensemble (Gaussian measure). Since then, a long list of further-reaching results about the entries of uniformly distributed random orthogonal or unitary  matrices  have been obtained. The most recent ones are   the previously cited paper of Diaconis {\it et al.},  the papers of Meckes and Chatterjee \cite{meckes08, meckessourav08}, the paper of Collins and Stolz \cite{collins-stolz08} and the paper of Jiang \cite{jiang06},  where the point of view is slightly different. In the present paper,   we give  a new, quite short, proof of the asymptotic normality of the linear combinations of the entries of uniformly distributed random unitary matrices, but we also extend   these investigations to the case where the distribution of the matrices is not the Haar measure but the heat kernel measure, with any initial distribution and any rescaling of the time.

 {\it Second order freeness:} A theory has been developed these last five years about Gaussian  fluctuations (called {\it second order limits}) of traces of large random matrices around their limits, the most emblematic articles in this theory   being \cite{mingo-nica04, mingo-speicher06, mingo-piotr-speicher07, mingo-piotr-collins-speicher07}. The results of this paper can be related to this theory, even though, technically speaking, we do not consider the powers of the matrices here\footnote{The reason is that the constant matrices we consider here, like $\sqrt{n}\times${\it (an elementary $n\times n$ matrix)}, have no bounded moments of order higher than two: our results are the best ones that one could obtain with such matrices.}. 

  {\it Brownian motion on the Lie algebra and It\^o map:}   The unitary Brownian motion    is  a continuous random process taking values on the unitary group, which has independent and stationary multiplicative increments. The most constructive way to define it is to
consider a  standard Brownian motion $(B_t)_{t\geq 0}$ on the tangent space of the unitary group at the identity matrix and to consider its image  by the {\it It\^o map} (whose inverse is sometimes called the {\it Cartan map}), i.e. to  wrap\footnote{\label{6309.11h55}For $G$ a matricial Lie group with tangent space $\mathfrak{g}$  at $I$, the ``wrapping" $w_\gamma$, on $G$, of a continuous and piecewise smooth  path $\gamma : [0,+\infty)\to \mathfrak{g}$ \st $\gamma(0)=0$ is defined  by $w_\gamma(0)=I$ and $w_\gamma'(t)=w_\gamma(t)\times \gamma'(t)$. It $B$ is a Brownian motion on $\mathfrak{g}$ and $(B_{n})_{n\geq 1}$ is a sequence of continuous, piecewise affine interpolations of $B$ with a step tending to zero as $n$ tends to infinity, then the sequence $w_{B_n}$ converges in \pro to a process which doesn't depend on the choice of the interpolations and which is a Brownian motion on $G$ 
\cite[Sect. VI.7]{iw81}, \cite[Eq. (35.6)]{rogers-williams2}, \cite{frizoberhauer}.} 
it around the unitary group: the  process $(U_t)_{t\geq 0}$  obtained  is a unitary Brownian motion   
starting at $I$. 
Our results give us an idea of the way the It\^o map alterates the Brownian feature of the entries of $(B_t)_{t\geq 0}$ at different scales of time.   Moreover, the question of the choice of a rescaling   of the time (depending on the dimension) raises interesting questions (see Remark \ref{27209.22h}).

\subsection*{Notation} For each $n\geq 1$, $\mathbb{U}_n$ shall denote the group of $n\times n$ unitary matrices. The identity matrix will always be denoted by $I$. For each complex matrix $M$,   $M^*$ will denote the adjoint of $M$. 
We shall call a standard complex Brownian motion a complex-valued process whose real and imaginary parts are independent   standard real Brownian motions divided by $\sqrt{2}$.
For all $k\geq 1$, the space of continuous functions from $[0,+\infty)$ to $\C^k$ will be denoted by $\mc{C}([0,+\infty),\C^k)$  and will be endowed with the topology of the uniform convergence on every compact interval.

\section{Statement of the results}
\subsection{Brief presentation of the Brownian motion on the unitary group}There are several ways to construct the Brownian motion on the unitary group\footnote{See \cite{hunt56, stroock-varadhan73,iw81,rogers-williams2}. A very concise and elementary definition  is also given in \cite{rains97}.}. For the one we choose here, all facts can easily be recovered  by the use of the matricial It\^o calculus, as exposed in Section \ref{27209.14h05}.

 Let $n$ be a positive integer and $\nu_0$ a \pro measure on the group of unitary $n\times n$ matrices. We shall call a {\it unitary Brownian motion with initial law $\nu_0$} any random process $(U_t)_{t\geq 0}$  with values on the space of $n\times n$ complex matrices \st $U_0$ is $\nu_0$-distributed and $(U_t)_{t\geq 0}$  is a strong solution of the stochastic differential equation \be\label{2220910h50} \ud U_t=i\ud H_tU_t-\ff{2}U_t\ud t,\ee
where $(H_{t})_{t\geq 0}$ is  a  Brownian motion\footnote{\label{1309.18h09}$(iH_t)_{t\geq 0}$ is in fact the skew-Hermitian Brownian motion that the process $(U_tU_{0}^*)_{t\geq0}$ wraps   around the unitary group, as explained in Footnote \ref{6309.11h55}.}
on the space of $n\times n$ Hermitian matrices endowed with the scalar product $\lan A,B\ran=n\Tr (AB)$.

It can  be proved that  for such a process $(U_t)_{t\geq 0}$,    for any $t_0\geq 0$, $U_{t_0}$ is almost surely  unitary and both processes $(U_{t_0+t}U_{t_0}^*)_{t\geq 0}$ and $(U_{t_0}^*U_{t_0+t})_{t\geq 0}$  are unitary Brownian motions starting at $I_n$ and independent of the  $\sigma$-algebra generated by $(U_s)_{0\leq s\leq t_0}$.

\begin{rmq}[Communicated by Thierry L\'evy]\label{27209.22h}{\rm There are other ways to scale the time for the Brownian motion on the unitary group. 
Our  scaling of the time is the one for which the three limit regimes correspond respectively to small values of $t$, finite values of $t$ and large values of $t$ and for which the limit non commutative distribution of $(U_t)_{t\geq 0}$ is the one of a free unitary Brownian motion.
It also has a heuristic  geometrical meaning: with this scaling, for any fixed $t$, the  distance\footnote{As explained in Footnotes \ref{6309.11h55} and \ref{1309.18h09}, $(U_t)_{t\geq 0}$ is the wrapping of $(iH_t)_{t\geq 0}$ on the unitary group, hence the distance between $U_0$ and $U_t$ has the same order as the one between $0$ and $iH_t$, which, by the Law of Large Numbers, has order $\sqrt{t}n$.} between $U_0$ and $U_t$  has the same order as  the diameter\footnote{By definition, the diameter of the group is the supremum, over pairs $U,V$ of unitary matrices, of the length of the shortest geodesic between $U$ and $V$. Here, it  is equal the length of the geodesic $t\in [0, n\pi]\mapsto \exp(itI_n/n)$  between  $I_n$ and  $-I_n$, i.e. to $n\pi$.}
 of the group.
 It means that for any fixed $t>0$,  large values of $n$,  $U_t$ is probably no longer too close to its departure point, while it also probably hasn't ``orbited" the unitary group  too many times.}\end{rmq}

\subsection{The three asymptotic regimes for the unitary Brownian motion starting at $I_n$}\label{2.04.09.11h30}

Let $(\al_n)_{n\geq 1}$ be a sequence of positive numbers with a limit $\al\in [0,+\infty]$.  
Let us fix a positive integer $k$ and let, for each $n\geq 1$, $A^n_{1}, \ldots, A^n_{k}$ be a family of non-random $n\times n$ matrices and  $(e^{-t/2}V^n_{t})_{t\geq 0}$ be a Brownian motion on $\mathbb{U}_{n}$  starting at $I$. Suppose that  there exists    complex matrices  $[a_l]_{ l=1}^k $, $[p_{l,l'}]_{ l,l'=1}^k$ and $[q_{l,l'}]_{ l,l'=1}^k$  \st for all $l,l'\in\{1, \ldots, k\}$, we have\beqy\label{cond.tr.25.01.09}{\ff{n}\Tr (A^n_{l})}&\ninf& a_{l},\\ \label{cond.tr.01.09}  \ff{n} \Tr (A^n_{l}A^n_{l'})& \ninf& p_{l,l'},\\
 \label{cond.norm.01.09}  {\ff{n}\Tr (A^n_{l}A^{n*}_{l'})} &\ninf &q_{l,l'} .\eeqy
 For each $n$, let us define, for $t\geq 0$, 
 $$X^n_{t}=\al_n^{-1/2}( \Tr[ A^n_{1}( V^n_{\log(\al_n t+1)}- I) ], \ldots, \Tr[ A^n_{k}(V^n_{\log(\al_n t+1)}- I) ]) .$$
Let $\mu$ be the \pro measure  on $\mc{C}([0,+\infty), \C^k)$ defined by the fact that any process $X_t=(X_{t,1},\ldots, X_{t,k})$ distributed according to $\mu$ is a  Gaussian centered process with independent increments  \st for all $t> 0$, for all $l,l'=1,\ldots, k$, \be\label{25309.11h12}\E(X_{t,l}\ovl{X_{t,l'}} )=q_{l,l'}t\,  , \; \E(X_{t,l}{X_{t,l'}} ) =\begin{cases}-p_{l,l'}t&\textrm{if $\al=0$,}\\
 -p_{l,l'}\f{\log(\al t+1)}{\al}+a_la_{l'}\f{\log^2(\al t+1)}{2\al}&\textrm{if $0<\al<+\infty$,}\\
 0&\textrm{if $\al=+\infty$.}\end{cases}\ee
 
  \begin{Th}\label{17209.0543}  As $n$ tends to infinity,  the distribution of the process $ (X^n_{t})_{t\geq 0} $ converges weakly to $\mu$.\end{Th}

\begin{rmq}\label{16.11.2010.0h50}{\rm Theorem \ref{17209.0543}  can easily be extended, 
  using standard topological arguments, to the case where for each $n$, the matrices $A_1^n, \ldots, A_k^n$ are random,  independent of  $(V^n_{t})_{t\geq 0}$ and the convergences of \eqref{cond.tr.25.01.09}, \eqref{cond.tr.01.09} and \eqref{cond.norm.01.09} hold in probability (with non-random limits).  It will be useful in the proofs of Theorem \ref{24209.07h44} and Theorem \ref{17209.0544}.}
\end{rmq}


Recall that a {\it principal submatrix} of a matrix is a matrix obtained by removing some columns, and the rows with the same indices.
\begin{cor}\label{26209.16h29} 
Let us fix $p\geq 1$  and let $(H_{t})$, $(S_t)$ be two independent standard Brownian motions on the  euclidian spaces of  $p\times p$ respectively Hermitian  and skew-Hermitian matrices endowed with the respective scalar products $\lan X,Y\ran=\Tr(XY)/2$, $\lan X,Y\ran=-\Tr(XY)/2$.   
Then, as $n$ tends to infinity,  the distribution of the $\C^{p\times p}$-valued process of the entries of any $p\times p$ principal submatrix  of $\sqrt{n/\al_n}(V^n_{\log(\al_n t+1)}- I)_{t\geq 0}$ converges   to the one of the random process $(H_{t-f_\al(t)}+S_{t+f_\al(t)})_{t\geq 0}$, where $$f_\al(t)=\begin{cases}t&\textrm{if $\al=0$},\\ \f{\log(\al t+1)}{\al}&\textrm{if $0<\al<+\infty$},\\
0&\textrm{if $\al=+\infty$.} 
\end{cases}$$\end{cor}


\begin{rmq}\label{20409.11h40}{\rm Note that when $\al=0$, the limit process is simply a standard Brownian motion on the space of $p\times p$ skew-Hermitian matrices, whereas, as $\al$ grows from zero to $+\infty$, the Hermitian part of the limit process keeps growing, and at last, when $\al=+\infty$, the limit process is a standard Brownian motion on the space of $p\times p$ complex matrices.  As said in the introduction, the existence of these three asymptotic regimes can be explained by the fact that 
the unitary Brownian motion is the ``wrapping", on the unitary group, of a Brownian motion on the tangent space  at $I$ (which is the space of skew-Hermitian matrices), and that as the time goes to infinity, its distribution tends to the Haar measure (for which, as stated by Corollary \ref{27309.23h50},  the upper-left corners are asymptotically distributed as standard complex Gaussian random matrices).}
\end{rmq}


\subsection{The particular case of  unitary Brownian motions with exchangeable rows and columns}
Let  $(\al_n)_{n\geq 1}$ be a sequence of positive numbers (no hypothesis is made on its convergence).  
Let us fix a positive integer $k$ and let, for each $n\geq 1$, $A^n_{1}, \ldots, A^n_{k}$ be a family of non-random $n\times n$ matrices and  $(e^{-t/2}V^n_{t})_{t\geq 0}$ be  a  Brownian motion on $\mathbb{U}_{n}$   \st $V^n_{0}$ is uniformly distributed on the group of matrices of permutations of $\{1,\ldots, n\}$. We suppose that  there exists  a   complex matrix $[q_{l,l'}]_{ l,l'=1}^k$   for all $l,l'\in \{1, \ldots, k\}$,  \beqy
 \label{20209.7h22} \ff{n}\Tr[A^n_{l}A^{n*}_{l'}]&\ninf &q_{l,l'},\\  \label{20209.7h23} \ff{n^2}\sharp \{(i,j)\ste \textrm{$((i,j)$-th entry of $A^n_{l})\neq 0$}\}& \ninf &0\eeqy
and that there exists a \pro measure $\mu_0$ on $\C^k$  \st \begin{align}
 \label{22.2.09.09h06}&& \al_n^{-1/2}(\Tr[A^n_{1}V^n_{0}],\ldots, \Tr[A^n_{k}V^n_{0}])\ninf \mu_0&&\textrm{(conv. in distribution).}
 \end{align}
 For each $n$, let us define, for $t\geq 0$,   
 $$X^n_{t}=\al_n^{-1/2}( \Tr(A^n_{1} V^n_{\log(\al_n t+1)}) , \ldots, \Tr(A^n_{k}V^n_{\log(\al_n t+1)}) ) .$$Let $P$ be an Hermitian matrix \st $P^2= [q_{l,l'}]_{l,l'=1}^k$, $(Z_{1},\ldots, Z_{k})$ be an independent family of standard complex Brownian motions and  $C$ be a $\mu_0$-distributed  random variable, independent of  the $Z_l$'s. Let us define $\mu$  to be the distribution, on $\mc{C}([0,+\infty), \C^k)$, of the process  $(C+(Z_{1,t},\ldots, Z_{k,t})P)_{t\geq 0}$.

   \begin{Th}\label{24209.07h44}      As $n$ tends to infinity, the distribution of $X_n$   converges weakly to $\mu$.\end{Th}
 
 \begin{rmq}\label{16.11.10.0h33}{\rm As Theorem \ref{17209.0543}, Theorem \ref{24209.07h44}  can   be extended 
    to the case where for each $n$, the matrices $A_1^n, \ldots, A_k^n$ are random,  independent of  $(V^n_{t})_{t\geq 0}$ and   the convergences of \eqref{20209.7h22} and  \eqref{20209.7h23} hold in probability  (with non-random limits). In several examples  given  below, the matrices $A_l^n$'s are actually random.}
\end{rmq}
 
 \begin{exs}{\rm Let us give a few examples of sequences $A^n$ which satisfy the hypotheses \eqref{20209.7h22}, \eqref{20209.7h23} and \eqref{22.2.09.09h06} (or their probabilistic versions mentioned in Remark \ref{16.11.10.0h33}). 
 
 a) Firstly, if $(\al_n)$ is bounded from below by a positive constant, if      \eqref{20209.7h22} holds and if for all $l$, $n^{-1}\sharp \{(i,j)\ste \textrm{$((i,j)$-th entry of $A^n_{l})\neq 0$}\}\longrightarrow 0$, then 
   \eqref{22.2.09.09h06}  holds   for $\mu_0$ the Dirac mass at zero (this can easily be deduced from Lemma \ref{23.2.09.07h35} bellow).
 
   b) Secondly, if, $\al_n$ tends to one and if, for all  $n$, the matrices $A^n_l$, $1\le l\le k$,  are    random real $n\times n$ matrices which satisfy  \eqref{20209.7h22} and \eqref{20209.7h23} for the convergence in probability, then \eqref{22.2.09.09h06} holds for $\mu_0$ the Gaussian measure with covariance matrix $[q_{l,l'}]_{l,l'=1}^k$. 
   This follows easily from \cite{schneller88}.
   As an example, if $k=1$ and if, for each $n$, the entries of $A^n_1$ are i.i.d. with distribution $\tau_n$ \st as $n$ tends to infinity, $$\tau_n(\{0\})\lto 1\qquad\textrm{ and }\qquad n(1- \tau_n(\{0\}))\int t^2\ud \tau_n(t)\lto 1,$$ then \eqref{20209.7h22}, \eqref{20209.7h23} and \eqref{22.2.09.09h06} hold for $q_{1,1}=1$ and $\mu_0$ the standard Gaussian law.
   
   c) Other examples can be found using \cite[Th. 5.1]{chen78}, where the laws $\mu_0$ are other infinitely divisible laws.}\end{exs}
 
 Both following corollaries  are direct applications of the previous theorem, the first one using implicitly the fact that  any entry of $V^n_0$ is null with \pro $1-n^{-1}$, and  the second one using implicitly   the fact that the distribution of the number of fixed points of a uniform random permutation of $\{1,\ldots, n\}$   converges weakly, as $n$ tends to infinity,    to the Poisson distribution with mean one \cite{ds94}.
 
\begin{cor}\label{26209.16h50}  Let $(\al_n)$ be a sequence of positive numbers. 
For any  $p,q\geq 1$, as $n$ tends to infinity,   the distribution of  any $p\times q$   submatrix of $$\lf(\sqrt{n/\al_n}V^n_{\log(\al_nt+1)}\ri)_{t\geq 0}$$ converges weakly  to the one of an independent family    of $pq$ standard    Brownian motions on the complex plane, i.e. a Brownian motion on the space of $p\times q$ complex matrices.
\end{cor}

\begin{cor} As $n$ tends to infinity,  the distribution of $(\Tr(V^n_{\log(t+1)}))_{t\geq 0}$ converges weakly  to the one of $(C+Z_t)_{t\geq 0}$, where $C$  is a Poisson random variable with mean one and $(Z_t)_{t\geq 0}$ is a standard complex Brownian motion, independent of $C$.
\end{cor}

\subsection{Application to the asymptotics of the uniform measure on the unitary group}
Since  
the Brownian motion on the unitary group distributed according to the Haar measure at time zero has a stationary distribution,  our results allow us to recover certain results of asymptotic normality of linear combinations of the entries of uniform random unitary matrices. 

The following theorem is not   new \cite{collins-stolz08,meckessourav08}. However, our method allows to give a very direct proof,  even under these very general hypotheses.
 \begin{Th}\label{17209.0544}  Let us fix $k\geq 1$ and let, for each $n\geq 1$, $A^n_{1},\ldots, A^n_{k}$
 be $n\times n$ non-random matrices and  $U_n$ be a random matrix with uniform distribution on the group of $n\times n$ unitary matrices. Suppose that   for all   $l,l'$, there is      $q_{l,l'}\in \C$ \st  $$\ff{n}\Tr(A^n_{l}A^{n*}_{l'})\ninf q_{l,l'}.$$ 
 Then as $n$ tends to infinity, the distribution of the random vector $(\Tr[A^n_{1}U_n],\ldots, \Tr[A^n_{k}U_n])$
 converges weakly to the one of   a Gaussian centered family $(Z_1,\ldots, Z_k)$ of complex random variables \st for all $l,l'=1,\ldots, k$, $\E(Z_lZ_{l'})=0$ and $ \E(Z_l\ovl{Z_{l'}})=q_{l,l'}$.
 \end{Th}
  
 The following corollary is immediate.
 \begin{cor}\label{27309.23h50}Let us fix $k\geq 1$ and let, for each $n$, $Z^n_{1},\ldots, Z^n_{k}$ be $k$ different entries of a random $n\times n$ matrix with uniform distribution on the unitary group. Then the joint distribution of $\sqrt{n}(Z^n_{1},\ldots, Z^n_{k})$ converges weakly, as $n$ tends to infinity, to the one of a family of independent standard complex Gaussian random variables. 
 \end{cor}

\section{Proofs}

\subsection{Preliminaries on matricial It\^o calculus}\label{27209.14h05}
a) Let $n$ be a positive integer. Let $(\mc{F}_t)_{t\geq 0}$ be a   filtration and 
$(H^n_{t})_{t\geq 0}$ be  an $(\mc{F}_t)_{t\geq 0}$-standard Brownian motion on the space of $n\times n$ Hermitian matrices endowed with the scalar product $\lan A,B\ran=n\Tr (AB)$, i.e. a process with values in the space of $n\times n$ Hermitian matrices \st the diagonal and upper diagonal 
entries of $(\sqrt{n}H^n_{t})_{t\geq 0}$ are independent random processes, the ones on the diagonal being standard real Brownian motions and the ones above the diagonal being standard complex Brownian motions. If one considers two 
 matrix-valued semimartingales  $X,Y$  \st $$\ud X_t =A_t(\ud H^n_{t})B_t+C_t\ud t, \quad\quad
\ud Y_t =D_t(\ud H^n_{t})E_t+F_t\ud t,$$ for some $(\mc{F}_t)_{t\geq 0}$-adapted matrix-valued processes $A,B, C,D,E, F$, then, by It\^o's formula, 
\beq\label{22.2.09.10h25}\ud (XY)_t&=&(\ud X_t)Y_t+X_t\ud Y_t+\ff{n}\Tr(B_tD_t)A_tE_t\ud t,\\
\label{22.2.09.10h26} \ud \lan \Tr(X),\Tr(Y)\ran_t&=&\ff{n}\Tr(B_tA_tE_tD_t)\ud t.\eeq
We shall use these formulas many times in the paper, without citing them every time.

b) With the same notation, let us now consider a deterministic $\mc{C}^1$ function $f$ with positive derivative \st $f(0)=0$.  Then for $X$ the process introduced above, the process $\tilde{X}_t:=X_{f(t)}$ satisfies $$\ud \tilde{X}_t =\sqrt{f'(t)}A_{f(t)}(\ud \tilde{H}^{n}_{t})B_{f(t)}+f'(t)C_{f(t)}\ud t,$$where $\tilde{H}^{n}$ is the $(\mc{F}_{f(t)})_{t\ge 0}$-Brownian motion defined by the formula $\tilde{H}^{n}_t=\int_0^t\ff{\sqrt{f'(s)}}\ud (B_{f(\cdot)})_s$.

\subsection{Proof of Theorem \ref{17209.0543}}Let us first state 
some matricial inequalities we shall often refer to in the following.  
Let $X,Y$ be two complex matrices and $G,H$ be two Hermitian nonnegative matrices. Then we have\beqy \label{matrix.10.2.09}|\Tr(XY)|&\leq &\sqrt{\Tr(XX^*)}\sqrt{\Tr (YY^*)},\\ \label{1329.matpos} \Tr(G^2)&\leq & (\Tr G)^2,\\\label{200309.12h31}|\Tr(GH)|&\leq &\Tr(G)\Tr(H).
\eeqy
Inequality \eqref{200309.12h31} follows from \eqref{matrix.10.2.09} and \eqref{1329.matpos}, which are obvious. 

   \begin{lem}\label{200109.0100} Let us fix $n\geq 3$, an  $n\times n$ matrix $A$  and
  a Brownian motion   $(e^{-t/2}V_{t})_{t\geq 0}$  on $\mathbb{U}_{n}$ starting at $I$. 
  Then there exists some real numbers $C_1,C_2,C_3 ,C_4$, independent of $t$,  whose absolute values are bounded by $100(\Tr AA^*)^2/n^2$   and \st  for all $t\geq 0$, \beqy
  \label{22309.09h} &&\E[\Tr(AV_tAV_t)]=\Tr(A^2)\cosh(t/n)-(\Tr A)^2\sinh(t/n),\\  
  \label{22309.09h01} && \E[|\Tr(AV_tAV_t)|^2]= |\Tr(A^2)|^2+\lf\{\f{|\Tr A|^4}{2n^4}+\f{C_1}{n^{3/2}}\ri\}\lf\{\cosh\lf(\f{2t}{n}\ri)-1\ri\}\\ 
  \nonumber&&-\lf\{\Re[(\Tr A)^2\ovl{\Tr(AA)}]+{n^{3/2}C_2}\ri\}\sinh\lf(\f{2t}{n}\ri) \\
  \nonumber&& +{nC_3}(e^t-1)+{C_4}(e^{2t}-1).
  \eeqy
   \end{lem}
  
    
    \begin{pr}Since the formulas we have to state are invariant under multiplication of $A$ by a scalar, we can suppose that $\Tr(AA^*)=n$. 
    
    Note that for  $H^n_{t}$ as in Section   \ref{27209.14h05}, $(V_t)$ is a strong  solution of    $\ud V_{t}=i(\ud H^n_{t}) V_{t}$. Hence by the matricial It\^o calculus,
     \beq
   \ud \Tr(AV_tAV_t)&=& 2i\Tr(V_tAV_tA\ud H_t^n)-{n}^{-1}\Tr(AV_t)\Tr(AV_t)\ud t,\\
    \ud \Tr(AV_t)\Tr(AV_t)&=& 2i\Tr(AV_t)\Tr(V_tA\ud H_t^n) -n^{-1}\Tr(AV_tAV_t)\ud t.
     \eeq
     It follows that for  $x(t)=\E[\Tr(AV_tAV_t)]$ and $y(t)=\E[\Tr(AV_t)\Tr(AV_t)]$, we have $$x'=-n^{-1}y\qquad\textrm{ and }\qquad y'=-n^{-1}x.$$ Equation \eqref{22309.09h} follows.

     Now, let us define, for   $C,D$ some $n\times n$ matrices, $u_{C,D}(t)=\E[\Tr(V_{t}CV_{t}^*D)]$ and $v_{C,D}(t)=\E[\Tr(V_{t}C)\Tr(V_{t}^*D)]$. By the matricial It\^o calculus
again, one has \be\label{uv200309}u_{C,D}(t)=\ff{n}(e^t-1)\Tr C \Tr D + \Tr(CD),\; v_{C,D}(t)=\ff{n}(e^t-1)\Tr(CD)+\Tr(C)\Tr(D).\ee
Let us now prove \eqref{22309.09h01}.
We introduce the functions \beq f(t)&=&\E(|\Tr(AV_tAV_t)|^2),\\ g(t)&=&\Re\{ \E[\Tr(AV_{t})\Tr(AV_{t})\Tr(A^*V_t^*A^*V_t^*)]\}, \\ h(t)&=&\E[|\Tr(AV_t)|^4].\eeq
    By the matricial It\^o calculus again (using   the hypothesis $\Tr(AA^*)=n$),
\beqy
\label{22309.12h30}n\times f'(t)&=&-{2}g(t)+{4e^t}u_{A^*A,AA^*}(t),\\
\nonumber n\times g'(t)&=&-f(t)-h(t)+4{e^t}\Re[v_{A,A^*AA^*}(t)],\\
\nonumber  n\times h'(t)&=&-2g(t)+4ne^tv_{A,A^*}(t).
\eeqy
It follows, by \eqref{uv200309}, that $g''(t)-\f{4}{n^2}g(t)=\f{8e^{2t}}{n}\kappa+\f{4e^{t}}{n}\theta$, for \beq \kappa&=&{n}^{-1}\Tr(AA^*AA^*)-1,\\
\theta&=&2-n^{-1}\Tr(AA^*AA^*)-n^{-1}\Tr(AAA^*A^*)-|\Tr A|^2+\Re\{\Tr A \Tr(A^*AA^*)\},\eeq hence $g(t)=\mu \cosh\lf(\f{2t}{n}\ri)+\nu \sinh\lf(\f{2t}{n}\ri)+\f{2n\kappa}{n^2-1}e^{2t}+\f{4n\theta}{n^2-4}e^t$, with 
\beq
\mu&=& \Re[\Tr A\Tr A\Tr(A^*A^*)]-\f{2n\kappa}{n^2-1}-\f{4n\theta}{n^2-4},\\
\nu&=&-\ff{2}|\Tr A^2|^2-\f{1}{2}|\Tr A|^4+{2}\Re[\Tr A\Tr(A^*AA^*)]-\f{2n^2\kappa}{n^2-1}-\f{2n^2\theta}{n^2-4}.
\eeq 
From \eqref{22309.12h30}, it follows that \bes\label{20309.14h11}f(t)=-|\Tr A^2|^2-\mu\sinh\lf(\f{2t}{n}\ri)-\nu\lf(\cosh\lf(\f{2t}{n}\ri)-1\ri)+w(t),\ees where
 $w(t)=-\f{2\kappa}{n^2-1}(e^{2t}-1)-\f{8\theta}{n^2-4}(e^t-1)+{2}(e^{2t}-1)+{4(n^{-1}\Tr(AA^*AA^*)-1)}(e^t-1)$.

Now, the conclusion follows from the fact that since $\Tr(AA^*)=n$, the inequalities \eqref{matrix.10.2.09}, \eqref{1329.matpos} and \eqref{200309.12h31} allow to prove that 
$|\Tr A|, |\Tr A^2| \leq n$, $\Tr (AA^*AA^*), |\Tr (AAA^*A^*)|\leq n^2$ and $|\Tr (A^*AA^*)|\leq {n^{3/2}}$.\end{pr} 

\begin{lem}\label{suite.20309}Let $ (a_n),(b_n), (c_n)$ be sequences of real numbers \st    $(a_n)$  tends to $+\infty$ and $(b_n)$ and $(c_n)$ are both bounded. Then we have  \be\label{19309.1} u_n:= \f{n^2b_n}{a_n}\lf(\cosh\lf(\f{\log a_n}{n}\ri)-1\ri)+ \f{nc_n}{a_n}\sinh\lf(\f{\log a_n}{n}\ri)
 \ninf 0.\ee \end{lem}
 
   \begin{pr}Let us define $$K=\max\lf\{\sup_{n\geq 1}|b_n|\, ,\;\sup_{n\geq 1} {|c_n|}\, ,\;\sup_{0<x\leq 1}\f{\cosh(x)-1}{x^2}\,,\;\sup_{0<x\leq 1}\f{\sinh(x)}{x}\ri\}.$$ Then \eqref{19309.1} follows from the fact that, since $u_n$ can also be written   
   $$u_n= \f{n^2b_n}{2a_n^{1-\ff{n}}}+ \f{n^2b_n}{2a_n^{1+\ff{n}}}- \f{n^2b_n}{a_n}
   + \f{nc_n}{2a_n^{1-\ff{n}}}- \f{nc_n}{2a_n^{1+\ff{n}}}\;,
   $$we have the upper-bound: $|u_n|\leq  5K\f{n^2}{e^{n-1}}\one_{a_n>e^{n}}+2K^2\f{\log^2 a_n +\log a_n}{a_n}\one_{a_n\leq e^{n}}.$
   \end{pr}

{\bf Proof of Theorem  \ref{17209.0543}.} 
For all $n$, $X^n$ is a $\C^k$-valued continuous centered martingale. To prove that its distribution tends to $\mu$, by Rebolledo's Theorem (see \cite{psv77} or \cite[Th. H.14]{alice-greg-ofer}), it suffices to prove that the bracket of $X^n$ converges pointwise, in $L^1$, to the one of a $\mu$-distributed process. Hence it suffices to fix     $\la_1, \ldots,\la_k \in \C$, to define the process $Y^n_t=\la_1X^n_{t,1}+\cdots+\la_n X^n_{t,k}$ and to prove that as $n$ tends to infinity,
 \be\label{22102010.19h51}\lan Y^n, \ovl{Y^n}\ran_t  
\ltoLtwo q t  \textrm{ and }\lan Y^n, {Y^n}\ran_t \ltoLtwo \begin{cases}-pt&\textrm{if $\al=0$,}\\
 -p\f{\log(\al t+1)}{\al}+a^2\f{\log^2(\al t+1)}{2\al}&\textrm{if $0<\al<+\infty$,}\\
 0&\textrm{if $\al=+\infty$.}\end{cases}\ee
for $q=\sum_{l,l'=1}^k\la_l\ovl{\la_{l'}}q_{l,l'}$, $p=\sum_{l,l'=1}^k\la_l{\la_{l'}}p_{l,l'}$ and $a=\sum_{l=1}^k\la_la_l$.

Let us define, for each $n$, $A^n=\sum_{l=1}^k\la_lA_l^n$. We have $$Y^n_t=
\al_n^{-1/2}\Tr[A^n(V^n_{\log(\al_n t+1)}-I)].$$Hence by Section   \ref{27209.14h05} b), $Y^n$ satisfies  $$\ud Y^n_t=\f{i}{\sqrt{\al_n t+1}}\Tr[V^n_{\log(\al_n t+1)}A^n\ud  {H}^n_t],$$ where ${H}^n_t$ is an Hermitian Brownian motion as introduced in Section \ref{27209.14h05}.
Thus, since $V^n_{\log(\al_n t+1)}V^{n*}_{\log(\al_n t+1)}=(\al_nt+1)I$, by the matricial It\^o calculus, we have  $$\ud \lan Y^n, \ovl{Y^n}\ran_t=\ff{n}\Tr(A^nA^{n*})\ud t,$$ so that  the first part of \eqref{22102010.19h51} follows directly from \eqref{cond.norm.01.09}. Let us now prove the second part. 
By the matricial It\^o calculus again, we have
$$\ud \lan Y^n, {Y^n}\ran_t=\f{-1}{n(\al_n t+1)}\Tr[A^nV^n_{\log(\al_n t+1)}A^{n}V^n_{\log(\al_n t+1)}]\ud t.$$Hence it suffice to prove that as $n$ tends to infinity, we have the convergence \be\label{22.4.10.23h09}\f{-1}{n(\al_n t+1)}\Tr[A^nV^n_{\log(\al_n t+1)}A^{n}V^n_{\log(\al_n t+1)}]\ltoLtwo \begin{cases}\f{-p}{\al t+1}+a^2\f{\log(\al t+1)}{\al t+1}&\textrm{if $0\le \al<+\infty$,}\\ 0&\textrm{if $\al=+\infty$,}\end{cases}\ee uniformly as $t$ varies in any compact subset of $[0, +\infty)$. This follows easily from Lemma \ref{200109.0100} (with Lemma \ref{suite.20309} in the case where $\al=+\infty$).\hfill$\square$

\subsection{Proof of corollary \ref{26209.16h29}}
By Theorem \ref{17209.0543} and formula \eqref{25309.11h12}, applied with matrices $A^n_{l}$ of the type $\sqrt{n}\times${\it (an elementary $n\times n$ matrix)}, the distribution  of any $p\times p$ principal submatrix  of $\sqrt{n/\al_n}(V^n_{\log(\al_n t+1)}- I)_{t\geq 0}$ converges  weakly   to the one of the random process $(M_t)_{t\geq 0}$ with  independent increments \st for all $p\times p$ complex matrices $X,Y$ and all $t\geq 0$,
$$\E[\Tr(M_tX)\ovl{\Tr(M_tY)}]=t\Tr(XY^*),\quad \E[\Tr(M_tX){\Tr(M_tY)}]=-\f{\log(\al t+1)}{\al}\Tr(XY),$$ where if $\al=0$ or $+\infty$, the second right-hand term has to be replaced by respectively $-t\Tr(XY)$ or $0$. Since a standard Brownian motion $(B_t)$ on an euclidian space $(E, \lan\cdot, \cdot\ran)$ satisfies, for all $u, v\in E$, $\E(\lan B_t, u\ran \lan B_t, v\ran)=\lan u,v\ran t$,  the result can easily be verified.\hfill$\square$

\subsection{Proof of Theorem \ref{24209.07h44}}

 \begin{lem}\label{23.2.09.07h35}
   Let us fix $n\geq 2$, let $S$ be the matrix of a uniform random permutation of $\{1,\ldots, n\}$ and let $A,B$ be  $n\times n$ non-random matrices. Then we have \beq\label{18.2.9.13h26} \E\{|\Tr (AS)|\}&\leq & \ff{n}\sqrt{C_A}\sqrt{\Tr(AA^*)},\\ \label{1829.13h27}
   \E\{|\Tr(ASBS)|\}&\leq & \f{n-1+\sqrt{C_AC_B}}{n(n-1)}\sqrt{\Tr(AA^*)}\sqrt{\Tr(BB^*)},\eeq where for each  matrix $X=[x_{i,j}]_{i,j=1}^n$,  $C_X$ denotes $\sharp\{(i,j)\ste x_{i,j}\neq 0\}$.
     \end{lem}
   
   \begin{pr}
Let us denote by respectively $[a_{i,j}], [b_{i,j}], [s_{i,j}]$ the entries of $A,B,S$.  
  By H\"older's inequality, we have $$ \E\{|\Tr (AS)|\}\leq \sum_{i,j} |a_{i,j}|\E(s_{i,j}) =\ff{n}\sum_{i,j}|a_{i,j}|\leq\ff{n}\sqrt{C_A}\sqrt{\Tr (AA^*)} .$$Moreover, we have  \beq
    \E\{|\Tr(ASBS)|\}
&\leq &\sum_{j,k}|a_{k,j}b_{k,j}|\E(s_{j,k}^2)+\sum_{\substack{i,j,k,l\\ j\neq l\textrm{ or } k\neq i}} |a_{i,j}b_{k,l}|\E(s_{j,k}s_{l,i})\\
&\leq &\ff{n}\sum_{j,k}|a_{k,j}b_{k,j}|+\ff{n(n-1)}\sum_{\substack{i,j,k,l}} |a_{i,j}b_{k,l}|,
\eeq 
and the conclusion follows from  H\"older's inequality again. 
\end{pr}

{\bf Proof of Theorem \ref{24209.07h44}.} We consider $C, Z_1,\ldots, Z_k,P$ as introduced above the statement of the theorem. For each   $t\geq 0$, let us define $K_t=(Z_{1,t},\ldots, Z_{k,t})P$.  It suffices to prove that as $n$ tends to infinity,  
the joint distribution of $(X^n_{0},(X^n_{t}-X^n_{0})_{t\geq 0})$ converges weakly to the one of $(C,(K_t)_{t\geq 0})$.

 First, Lemma \ref{23.2.09.07h35} and the hypothesis  \eqref{20209.7h22} and \eqref{20209.7h23} allow us to claim that as $n$ tends to infinity,  $\ff{n}\Tr(A_{n,l}V_{n,0})$ and $ \ff{n}\Tr(A_{n,l}V_{n,0}A_{n,l'}V_{n,0})$ both converge to zero in probability. 
 Moreover, by a standard topological argument, one can suppose that $(\al_n)$ admits a limit  $\al\in [0,+\infty]$. 
 
 Now, note that for all $n$, $(e^{-t/2}V_{0}^{n*}V^n_{t})_{t\geq 0}$ is a unitary Brownian motion starting at $I$, independent of $V^n_0$. 
 By Theorem \ref{17209.0543} and Remark \ref{16.11.2010.0h50}, it implies that the joint distribution of $(X^n_{0},(X^n_{t}-X^n_{0})_{t\geq 0})$ converges weakly to the one of $(C,(K_t)_{t\geq 0})$, which closes the proof of the theorem.\hfill$\square$

\subsection{Proof of Theorem \ref{17209.0544}}
\begin{lem}\label{16209.10h}Let $U$ be  Haar-distributed on $\mathbb{U}_n$  and $A$ be an $n\times n$ matrix, with $n\ge 3$. Then
 \beqy\label{1520912h00} \E\{|\Tr (AU)|^2\}&=&\Tr(AA^*)/n,\\ \label{1520912h01}\E\{|\Tr(AUAU)|^2\}&\leq &{100}(\Tr(AA^*))^2/n^2.\eeqy
 \end{lem}
 
 \begin{pr} Set $U=[u_{i,j}]_{i,j=1}^n$.
One can write $A=VA'W$, with $V,W$ unitary matrices and $A'$ a diagonal matrix whose diagonal entries $a_1,\ldots, a_n$   are the eigenvalues of $\sqrt{AA^*}$. Since the law of $U$ is invariant under the left and right actions of the unitary group, one can suppose that $A=A'$.  These invariances of the law of $U$ also imply that for all $i,j$, $\E(u_{i,i}\overline{u_{j,j}})=\delta_i^j/n$. Equation \eqref{1520912h00} follows.
Equation  \eqref{1520912h01} follows from \eqref{22309.09h01} and the fact that the Haar measure on the unitary group is the limit of the distribution of $e^{-t/2}V_t$ as $t$ tends to infinity.
 \end{pr}
 
{\bf Proof of Theorem \ref{17209.0544}.}
{\it Step I.} Firstly,  by  Lemma \ref{16209.10h}, for all $l,l'$,  $n^{-1}\Tr(A_{l}^nU_n)$ and $n^{-1}\Tr(A_{l}^nU_nA_{l'}^nU_n)$ both tend in \pro  to zero as $n$ tends to infinity. 

  {\it Step II.}  Following \cite[Th.  D8]{DemboZeitouni},  we fix a bounded real function $f$ on $\C^k$ which is $1$-Lipschitz  for the canonical hermitian norm $||\cdot ||$ on $\C^k$, and we shall prove that \be\label{1620910h23}\E\{f(\Tr[A^n_{1}U_n],\ldots, \Tr[A^n_{k}U_n])\}\ninf \E\{f(Z_1,\ldots, Z_k)\}.\ee 
  Let us fix $\eps >0$. 
  
a)  Let, for each $t\geq 0$, $(Z_{1,t},\ldots, Z_{k,t})$ be a  Gaussian family of  centered complex random variables \st  for all $l,l'$, $\E[Z_{l,t}Z_{l',t}]=0$ and   $\E[Z_{l,t}\ovl{Z_{l',t}}]=q_{l,l'}(1-e^{-t})$ (such a family exists because the matrix $[q_{l,l'}]_{l,l'}$ is   nonnegative). 
The distribution of $(Z_{1,t},\ldots, Z_{k,t})$ tends to the one of $(Z_1,\ldots, Z_k)$ as $t$ tends to infinity. Hence there is $t_0>0$ \st \beqy \label{160209.23h} |\E\{f(Z_{1},\ldots, Z_{k})\}-\E\{f(Z_{1,t_0},\ldots, Z_{k,t_0})\}|&\leq &\eps,\\  \label{160209.23h.primabord} 
e^{-\f{t_0}{2}}\sup_{n\geq 1}\{n^{-1}\Tr[A_{1}^nA_{1}^{n*}+\cdots +A^n_{k}A_{k}^{n*}]\}^{{1}/{2}}&\leq &\eps.
\eeqy
\indent b)  For each $n$, up to an extension of the \pro space where   $U_n$ is defined, one can suppose that there exists a unitary Brownian motion  $(U^n_{t})_{t\geq 0}$, starting at $I$, independent of $U_n$. Let us define, for each $n$,\beq X^n&=& (\Tr[A^n_{1}U_nU^n_{t_0}],\ldots, \Tr[A^n_{k}U_nU^n_{t_0}])\\
 Y^n&=& e^{-\f{t_0}{2}}(\Tr[A^n_{1}U_n],\ldots, \Tr[A^n_{k}U_n])\\
 D^n&=&X^n- Y^n=e^{-\f{t_0}{2}}(\Tr[A^n_{1}U_n(e^{\f{t_0}{2}}U^n_{t_0}-I)],\ldots,\Tr[A^n_{k}U_n(e^{\f{t_0}{2}}U^n_{t_0}-I)]) \eeq
 By Step I and the randomized version of Theorem \ref{17209.0543} stated in Remark \ref{16.11.2010.0h50}, as $n$ tends to infinity, the distribution of  $D^n$ converges to  the one of 
 $(Z_{1,t_0},\ldots, Z_{k,t_0})$. It follows that for $n$ large enough,\be\label{17209.05h54}|\E\{f(Z_{1,t_0},\ldots, Z_{k,t_0})\}-\E\{f(D^n)\}|\leq \eps.
 \ee
 
c) At last, since $f$ is $1$-Lipschitz   for $||\cdot ||$, for all $n$, we have \beq|\E\{f(D^n)\}-\E\{f(X^n)\}|&\leq &\E\{||D^n -X^n||\}\\
 &\leq & (e^{-t_0}\E\{ |\Tr[A^n_{1}U_n]|^2+\cdots+|\Tr[A^n_{k}U_n]|^2\})^{1/2}\\
 &\leq & e^{-t_0/2}(n^{-1}\Tr[A^n_{1}A^{n*}_{1}]+\cdots+n^{-1}\Tr[A^n_{k}A^{n*}_{k}])^{1/2},
 \eeq
 the last inequality following from \eqref{1520912h00}.
 By \eqref{160209.23h.primabord},  
  it allows us to claim that for all $n$, \be\label{17209.06h14}|\E\{f(D^n)\}-\E\{f(X^n)\}|\leq\eps.\ee
   
d) To conclude, note that by the right invariance of the law of $U_n$, $X^n$ and the random vector   $(\Tr[A^n_{1}U_n],$ \ldots,$\Tr[A^n_{k}U_n])$  have the same distribution. By  \eqref{160209.23h}, \eqref{17209.05h54} and \eqref{17209.06h14}, it follows that for 
$n$ large enough, $$| \E\{f(Z_1,\ldots, Z_k)\}-\E\{f(\Tr[A^n_{1}U_n],\ldots, \Tr[A^n_{k}U_n])\}|\leq 3\eps.$$ It closes the proof of the theorem.\hfill$\square$

\end{document}